\documentclass[centertags,12pt]{amsart}
\usepackage{latexsym}
\usepackage{amssymb}
\numberwithin{equation}{section}
\makeatletter
\renewcommand{\subsection}{\@startsection
{subsection}{2}{0mm}{\baselineskip}{-0.25cm}
{\normalfont\normalsize\bf}}
\makeatother
\newtheorem*{theorem*}{Theorem}
\newtheorem{theorem}{Theorem}[section]
\newtheorem{proposition}[theorem]{Proposition}
\newtheorem{corollary}[theorem]{Corollary}
\newtheorem{lemma}[theorem]{Lemma}
\theoremstyle{definition}
\newtheorem{remark}[theorem]{Remark}
\def\P{\mathbb P}

\def\cC{\mathcal C}
\def\cD{\mathcal D}
\def\cF{\mathcal F}
\def\cH{\mathcal H}

\def\cX{\mathcal X}
\def\cY{\mathcal Y}
\def\cU{\mathcal U}
\def\cV{\mathcal V}
\def\cT{\mathcal T}

\def\cZ{\mathcal Z}
\def\sq{\sqrt{q}}
\def\sqq{\sqrt{q} + 1}
\def\fq{\mathbb F_q}

\def\div{{\rm div}}
\def\dim{{\rm dim}}
\def\deg{{\rm deg}}
\def\det{{\rm det}}
\def\ord{{\rm ord}}
\def\frx{\Phi_{\mathcal X}}
\def\qa{\textstyle\frac 14}
\def\ha{\textstyle\frac 12}
\def\ei{\textstyle\frac 18}
\begin{document}
\author[Cossidente]{A. Cossidente}
\author[Hirschfeld]{J.W.P. Hirschfeld}
\author[Korchm\'aros]{G. Korchm\'aros}
\author[Torres]{F. Torres}\thanks{1991 Math. Subj. Class.:
Primary 11G, Secondary 14G}
\title[Plane maximal curves]{On plane maximal curves}
\address{Dipartimento de Matematica, Universit\`a della Basilicata,
Potenza, 85100, Italy}
\email{cossidente@unibas.it}
\address{School of Mathematical Sciences, University of Sussex, Brighton
BN1 9QH, United Kingdom}
\email{jwph@sussex.ac.uk}
\address{Dipartimento de Matematica, Universit\`a della Basilicata,
Potenza, 85100, Italy}
\email{korchmaros@unibas.it}
\address{IMECC-UNICAMP, Cx. P. 6065, Campinas, 13083-970-SP, Brazil}
\email{ftorres@ime.unicamp.br}
\begin{abstract}
The number $N$ of rational points on an algebraic curve of genus $g$ over
a finite field $\fq$ satisfies the Hasse-Weil bound $N \leq q + 1 +
2g\sq$.  A curve that attains this bound is called maximal. With
 $g_0 =\ha (q -\sq)$ and $g_1=\qa (\sq - 1)^2$, it is known that maximal
curves have $g= g_0$ or $g \leq g_1$. Maximal curves with $g=g_0$ or
$g_1$ have been characterized up to isomorphism. A natural genus to be
studied  is 
$$
g_2 =\ei (\sq - 1)(\sq - 3)\, ,
$$ 
and for this genus there are two non-isomorphic maximal curves known when
$\sq\equiv 3\pmod{4}$. Here, a maximal curve with genus $g_2$ and a
non-singular plane model is characterized as a Fermat curve of degree
$\ha(\sqq)$. 
\end{abstract}
\maketitle
\section{Introduction}\label{s1}
For a non-singular model of a projective, geometrically irreducible,
algebraic curve $\cX$ defined over a finite field $\fq$ with $q$ elements,
the number $N$ of its $\fq$-rational points satisfies the Hasse-Weil
bound, namely (see \cite{weil}, \cite[\S V.2]{sti}) 
$$
|N-(q+1)|\le 2g\sqrt{q}\, .
$$
If $\cX$ is plane of degree $d$, then this bound implies that
\begin{equation}\label{eq1}
|N-(q+1)|\le (d-1)(d-2)\sqrt{q}\, .
\end{equation}
These bounds are important for applications in Coding Theory (see, for
example, \cite{sti}) and in Finite Geometry (see \cite[Ch. 10]{h}). In
these subjects one is often interested in curves with {\em many}
$\fq$-rational points and, in particular, {\em maximal curves}, that is,
curves where $N$ reaches the upper Hasse-Weil bound. 

The approach of St\"ohr and Voloch \cite{sv} to the Hasse-Weil bound
shows that an upper bound for $N$ can be obtained via $\fq$-linear
series. This upper bound depends not only on $q$ and $g$, as does the
Hasse-Weil bound, but also on the dimension and the degree of the linear
series.

In \cite{hk1} an upper bound for $N$ was found in the case that $\cX$ is a
plane curve. It turns out that this bound is better than the upper bound
from (\ref{eq1}) under certain conditions on $d$ and $q$. The bound in
\cite{hk1} is not symmetrical in the different types of branches. Two
types of branches are distinguished, both centred at $\fq$-rational points
of $\cX$: (a) the branches of order $r$ and class $r$; (b) the branches of
order $r$ and class different from $r$. Let $M_q$ and $M_q'$ be the number
of branches of type (a) and (b)  respectively, each counted $r$ times. If
$d$ and $q$ satisfy certain restrictions, then
\begin{equation}\label{eq1.2'} 
2M_q+M_q'\le d(q-\sq+1)\, , 
\end{equation} 
and equality holds if and only if $\cX$ is a non-singular plane maximal
curve over $\fq$ of degree $\ha(\sq+1)$. This result is the starting point
of our research. 

An example of a curve attained the equality in (\ref{eq1.2'}) is provided
by the Fermat curve $\cF$ (see \S\ref{s3}) with equation, in homogeneous
coordinates $(U,V,W)$,
\begin{equation}\label{eq1.2}
U^{(\sq+1)/2}+V^{(\sq+1)/2}+W^{(\sq+1)/2}=0\, .
\end{equation}
The main result of the paper is to show the following converse 
(see \S\ref{s5}).
\begin{theorem}\label{thm1} 
If $\cX$ is a non-singular plane maximal curve over $\fq$  
of degree $\ha(\sq+1),$ then it is $\fq$-isomorphic to $\cF$.
\end{theorem}
This result is connected to recent investigations on the genus of maximal
curves \cite{ft}, \cite{fgt}, \cite{ft1}. The genus $g$ of a maximal curve
$\cX$ over $\fq$ is at most $\ha\sq(\sq-1)$ \cite{ihara}, \cite[\S
V.2]{sti} with equality holding if and only if $\cX$ is $\fq$-isomorphic
to the Hermitian curve with equation
$$
u^{\sq+1}+v^{\sq+1}+w^{\sq+1}=0\, ,
$$
\cite{r-sti}. In \cite{ft} it  was observed that
$$ 
g\le \qa(\sq-1)^2\quad \mbox{ if }\quad g< \ha\sq(\sq-1)\, ,
$$ 
a result conjectured in \cite{sti-x}. Also, if $q$ is odd and 
$$
\qa(\sq-1)(\sq-2) < g \le \qa(\sq-1)^2\, ,
$$ 
then $g=\qa(\sq-1)^2$ and $\cX$ is
$\fq$-isomorphic to the non-singular model of the curve with affine
equation 
$$
y^q+y=x^{(\sq+1)/2}\, ,
$$
\cite[Thm. 3.1]{fgt}, \cite[Prop. 2.5]{ft1}. In general, the
situation for either $q$ odd and $g\le \qa(\sq-1)(\sq-2)$ or $q$ even and
$g\le \qa\sq(\sq-2)$ is unknown. In the latter case, an example
where equality holds is provided by the non-singular model of the curve 
with affine equation
$$
\sum_{i=1}^{t}y^{\sq/2^i}= x^{q+1}\, ,\quad \sq=2^t\, ,
$$ 
and it seems that this example may be the only one up to $\fq$-isomorphism 
\cite{at}.

In \cite[\S2]{fgt} the maximal curves obtained from the affine equation 
$$
y^{\sq}+y=x^m\, ,
$$ 
where $m$ is a divisor of $(\sq+1)$, are characterized by means of
Weierstrass semigroups at an $\fq$-rational point; the genus of these
curves is $g=\ha(\sq-1)(m-1)$. If $m=\qa(\sq+1)$ and $\sq\equiv
3\pmod{4}$, we find two curves of genus $\ei(\sq-1)(\sq-3)$, namely the
curve with affine equation
$$
y^{\sq}+y=x^{(\sq+1)/4}
$$ 
and the curve $\cX$ of our main result. It turns out that these curves are
not $\bar\fq$-isomorphic (see Remark \ref{rem4.1}(ii)). As far as we know,
this is the first example of two maximal curves of a given genus that are
not $\fq$-isomorphic for infinitely many values of $q$. 

As in \cite{hk}, \cite{hk1}, \cite{ft}, \cite{fgt}, \cite{ft1}, the key
tool used to carry out the research here is the approach of St\"ohr and
Voloch \cite{sv} to the Hasse-Weil bound applied to suitable $\fq$-linear
series on the curve. 

\noindent
{\bf Convention.} From now on, the word {\em curve} means a projective,
geometrically irreducible, non-singular, algebraic curve.
\section{Background}\label{s2}
In this section we summarize background material concerning
Weierstrass points and Frobenius orders from \cite[\S\S1--2]{sv}.

Let $\cX$ be a curve of genus $g$ defined over $\bar\fq$ equipped with the
action of the Frobenius morphism $\frx$ over $\fq$. Let $\cD$ be a $g^r_d$
on $\cX$ and suppose that it is defined over $\fq$. Then associated to
$\cD$ there exist two divisors on $\cX$, namely the {\em ramification
divisor}, denoted by $R=R^{\cD}$, and the {\em $\fq$-Frobenius divisor},
denoted by $S=S^\cD=S^{(\cD,q)}$. Both divisors describe the
geometrical and arithmetical properties of $\cX$; in particular, 
the divisor $S$ provides information on
the number $\# \cX(\fq)$ of $\fq$-rational points of $\cX$.

For $P\in \cX$, let $j_i(P)$ be the $i$th $(\cD,P)$-order,
$\epsilon_i=\epsilon_i^\cD$ be the $i$-th $\cD$-order ($i=0,\ldots,r$), and
$\nu_i=\nu_i^{(\cD,q)}$ be the $i$-th $\fq$-Frobenius order of $\cD$
($i=0,\ldots,r-1$). Then the following properties hold:
\begin{enumerate}
\item\label{1} $\deg(R)=(2g-2)\sum_{i=0}^{r}\epsilon_i + (r+1)d$;
\item\label{2} $j_i(P)\ge \epsilon_i$ for each $i$ and each $P$;
\item\label{3} $v_P(R)\ge \sum_{i=0}^{r} (j_i(P)-\epsilon_i)$ and equality
holds if and only if $\det(\binom{j_i(P)}{\epsilon_j}) \not\equiv 0\pmod{p}$;
\item\label{4} $(\nu_i)$ is a subsequence of $(\epsilon_i)$;
\item\label{5} $\deg(S)= (2g-2)\sum_{i=0}^{r-1}\nu_i + (q+r)d$;
\item\label{6} For each $i$ and for each $P\in \cX(\fq)$, $\nu_i\le
j_{i+1}(P)-j_1(P)$;
\item\label{7} For each $P\in \cX(\fq)$, $v_P(S)\ge
\sum_{i=0}^{r-1}(j_{i+1}(P)-\nu_i)$ and equality holds if and only if 
$\det(\binom{j_{i+1}(P)}{\nu_j})\not\equiv 0 \pmod{p}$.
\end{enumerate}
Therefore, if $P\in \cX(\fq)$, properties (6) and (7) imply
\begin{enumerate}
\item[8.] $v_P(S)\ge rj_1(P)$.
\end{enumerate}
Consequently, from (5) and (8), we obtain the main result of
\cite{sv}, namely
\begin{enumerate}
\item[9.] $\# \cX(\fq)\le \deg(S)/r$.
\end{enumerate}
\section{Plane maximal curves of degree $(\protect\sq+1)/2$} \label{s3}
Throughout this section we use the following notation:
\begin{enumerate}
\item[(a)]  $\Sigma_1$ is the linear series on a plane
curve over $\fq$ obtained from lines  of $\P^2(\fq)$,
and  $\Sigma_2$ is the series obtained from conics;
\item[(b)] for
$i=1,2$, the divisor $R_i$ is the ramification divisor and  $S_i$ 
is the $\fq$-Frobenius divisor associated to $\Sigma_i$.;
\item[(c)]
$j^i_n(P)$ is the $n$-th $(\Sigma_i,P)$-order;
\item[(d)]
$\epsilon^i_n=\epsilon^{\Sigma_i}_n$ and $\nu^i_n=\nu^{(\Sigma_i,q)}_n$;
\item[(e)]  $p={\rm char}(\fq)$.
\end{enumerate}
\begin{lemma}\label{lemma3.1} Let $\cX$ be a plane non-singular curve
over $\fq$ of degree $d$. 
If $d\not\equiv 1\pmod{p},$ then $\cX$ is
classical for $\Sigma_1$.
\end{lemma}
\begin{proof} See \cite[Corollary 2.2]{pardini} for $p>2$, and
\cite[Corollary 2.4]{homma} for $p\ge 2$ .
\end{proof}
\begin{corollary}\label{cor3.1} Let $\cX$ be a plane non-singular
maximal curve over $\fq$ of degree $d$ with $d\not\equiv 0\pmod{p}$ and
$2<d< (\sq+1)^2/3$. Then there exists
$P_0\in \cX(\fq)$ whose $(\Sigma_1,P_0)$-orders are $0, 1, 2$.
\end{corollary}
\begin{proof} Suppose that $j^1_2(P)>2$ for each $P\in \cX(\fq)$. Then
by  \S\ref{s2}(\ref{3}) and the previous lemma we would have $v_P(R_1)\ge
1$ for such points $P$.
Consequently, by \S\ref{s2}(\ref{1}) and the maximality of $\cX$ it
follows that
$$
{\rm deg}(R_1)=3(2g-2)+3d\ge \# \cX (\fq)= (\sq+1)^2+\sq(2g-2)
$$
so that
$$
0\ge (\sq+1)(\sq+1-\frac{3d}{\sq+1})+(2g-2)(\sq-3)\, ,
$$
a contradiction.
\end{proof}

Throughout the remainder of the paper, let $\cX$ be a plane non-singular maximal curve of degree $d$.
We have the following relation between $(\Sigma_1,P)$-orders and
$(\Sigma_2,P)$-orders for $P\in \cX$.
\begin{remark}[GV, p. 464]\label{rem3.1} For $P\in \cX$, the following
set
$$
\{j^1_1(P), j^1_2(P), 2j^1_1(P), j^1_1(P)+j^1_2(P), 2j^1_2(P)\}
$$
is contained in the set of $(\Sigma_2,P)$-orders.
\end{remark}

Now suppose that $d$ satisfies the hypotheses in Corollary \ref{cor3.1}
and let $P_0\in \cX(\fq)$ be as in this corollary. Then, by Remark
\ref{rem3.1} and the fact that $\dim(\Sigma_2)=5$, the $(\Sigma_2,P_0)$-orders
are $0,1,2,3,4$ and $j:=j^2_5(P_0)$ with $5\le j\le \sq$. Therefore, by
\S\ref{s2}(\ref{2}),(\ref{6}),(\ref{4}),
\begin{enumerate}
\item[(a)] the $\Sigma_2$-orders are
$0,1,2,3,4$ and $\epsilon:=\epsilon^2_5$ with $5\le\epsilon\le j$;
\item[(b)]  the $\fq$-Frobenius orders are
$0,1,2,3$ and $\nu:=\nu^2_4$ with $\nu\in\{4, \epsilon\}$.
\end{enumerate}
\begin{corollary}\label{cor3.2} Let $\cX$ be a plane non-singular maximal
curve over $\fq$ of degree\\
$d=\ha(\sq+1)$. If $\sq\ge 11,$ then
\begin{enumerate}
\item[\rm(1)] the $\Sigma_2$-orders are $0, 1, 2, 3, 4,\sq;$
\item[\rm(2)] the $\fq$-Frobenius orders of $\Sigma_2$ are $0, 1, 2, 3,\sq$.
\end{enumerate}
\end{corollary}
\begin{proof} The curve $\cX$ satisfies the hypotheses in Corollary
\ref{cor3.1}. So, with the above notation, we have to show that
$\epsilon=\nu=\sq$.

(a) First it is shown that  $\nu=\epsilon$.

 We have already seen that
$\nu\in\{4,\epsilon\}$. {}From
\S\ref{s2}(\ref{5}),(8) and the maximality of $\cX$ we have that
$$
\deg(S_2)=(6+\nu)(2g-2)+(q+5)(\sq+1)\ge 5\#\cX(\fq)=5(\sq+1)^2 +5\sq(2g-2),
$$
so that
\begin{equation}\label{eq3.1}
(\sq-5)(\sq-6-\nu)\le 0\, .
\end{equation}
Then if $\nu=4$, we would have $\sq\le 10$, a contradiction.

(b)
Now, $p$ divides $\epsilon$ (see \cite[Corollary
3]{garcia-homma}).
From \S\ref{s2}(\ref{6}) and (a),
$$
\nu=\epsilon\le j_5(P_0)-j_1(P_0)\le \sq\, .
$$
Therefore,
from (\ref{eq3.1}), the fact that $\sq>5$ and (a),
$$
\epsilon\in\{\sq-6,\sq-5,\sq-4,\sq-3,\sq-2\,\sq-1,\sq\}\, .
$$
Since $p > 2$ and $p$ divides $\epsilon$, the possibilities are reduced to
the following:
$$
\epsilon\in\{\sq-6,\sq-5,\sq-3,\sq\}\, .
$$

If $\epsilon=\sq-6$, then $p=3$ and by the $p$-adic criterion
\cite[Corollary 1.9]{sv} $\epsilon=6$ and so $\sq=12$, a contradiction.

If $\epsilon=\sq-5$, then $p=5$. Since $\binom{\sq-5}{5}\not\equiv
0\pmod{5}$, by the $p$-adic criterion we would have that 5 is
also a $\Sigma_2$-order, a contradiction.

If $\epsilon=\sq-3$, then $p=3$ and so $\sq=9$, which is eliminated by the
hypothesis that $\sq\ge 11$.

Hence $\epsilon=\sq$, which completes the proof.
\end{proof}

Now the main result of this section can be stated. We recall that a
maximal curve $\cX$ over $\fq$ is equipped with the $\fq$-linear series
$\cD_\cX:=|(\sq+1)P_0|$, $P_0\in \cX(\fq)$, which is independent of
$P_0$ and provides a lot of information about the curve (see
\cite[\S1]{fgt}).
\begin{theorem}\label{thm3.1} Let $\cX$ be a plane maximal curve over
$\fq$ of degree $\ha(\sq+1)$. Suppose that $\sq\ge 11$. Then the linear
series $\cD_\cX$ is the linear series
$\Sigma_2$ cut out by conics.
\end{theorem}
\begin{proof} First it is shown that, for $P\in \cX(\fq)$, the
intersection divisor of the
osculating conic $\cC^{(2)}_P$ and $\cX$ satisfies
\begin{equation}\label{eq3.2}
\cC^{(2)}_P. \cX =(\sq+1)P.
\end{equation}
To show this, let $P\in X(\fq)$; then, by Corollary \ref{cor3.2}(1) and
\S\ref{s2}(\ref{6}), we have that $\nu=\sq\le j_5(P)-j_1(P)\le \sq$
(recall that $\deg(\Sigma_2)=\sq+1$). Consequently $j^2_5(P)=\sq+1$ and
so (\ref{eq3.2}) follows.

This implies that $\Sigma_2\subseteq \cD$. Then to show the equality it is
enough to show that $n+1:=\dim(\cD)\le 5$. To see this we use
Castelnuovo's genus bound for curves in projective spaces as given in
\cite[p.34]{fgt}: the genus $g$ of $\cX$
satisfies
$$
2g\le \begin{cases}
(2\sq-n)^2/(4n) & \text{if $n$ is even}, \\
((2\sq-n)^2-1)/(4n) & \text{if $n$ is odd}.
\end{cases}
$$
Suppose that $n+1\ge 6$. Then, since $2g-2=(\sq-1)(\sq-3)/4$, we would
have
$$
(\sq-1)(\sq-3)/4\le ((2\sq-5)^2-1)/20=(\sq-3)(\sq-2)/5\, ,
$$
a contradiction. This finishes the proof.
\end{proof}

Next we compute the $(\Sigma_1,P)$-orders for $P\in \cX$.
\begin{lemma}\label{lemma3.2} Let $\cX$ be a plane maximal curve over
$\fq$ of degree $\ha(\sq+1)$ and let $P\in \cX$.
\begin{enumerate}
\item[\rm(1)] Two types of $\fq$-rational points of $\cX$ are distinguished$:$

{\rm(a)} regular points$,$ that is$,$ points whose $(\Sigma_1,P)$-orders
are $0, 1, 2,$ so that\\ $v_P(R_1)=0;$

{\rm(b)} inflexion points$,$ that is$,$ points whose 
$(\Sigma_1,P)$-orders are $0, 1,\ha(\sq+1),$\\ so that $v_P(R_1)=(\sq-3)/2$.
\item[\rm(2)] If $P\not\in \cX(\fq),$ then the $(\Sigma_1,P)$-orders
are
$ 0, 1,2,$ so that $v_P(R_1)=0$.
\end{enumerate}
\end{lemma}
\begin{proof} For each $P\in\cX$ we have that $j^1_1(P)=1$ because $\cX$
is non-singular. So we just need to compute $j(P):=j^1_2(P)$.

We know that $\cD_\cX=\Sigma_2=2\Sigma_1$, $\dim(\Sigma_2)=5$, and that
$j^2_5(P)=\sq+1$ provided that $P\in \cX$ (see proof of Theorem
\ref{thm3.1}). In addition, by \cite[Thm. 1.4(ii)]{fgt}, $j^2_5(P)=\sq$
for $P\not\in\cX(\fq)$.

Suppose that $j(P)>2$. Then from Remark \ref{rem3.1} we must have
$j^2_5(P)=2j(P)$. Since $\sq$ is odd, this is the case if and only if $2j(P)=\sq+1$
and $P\in \cX(\fq)$, because of the above computations.

The computations for $v_P(R_1)$ follow from \S\ref{s2}(\ref{3}).
\end{proof}

Let
\begin{align*}
M_q=M_q(\cX):= & \# \{P\in \cX(\fq): j^1_2(P)=2\}\, ,\\
\intertext{and}
M_q'=M_q(\cX)':= & \# \{P\in \cX(\fq): j^1_2(P)=\ha(\sq+1)\}\, .
\end{align*}
\begin{theorem}\label{thm3.2} Let $\cX$ be a plane maximal curve over
$\fq$ of degree $\ha(\sq+1)$. Suppose that $\sq\ge 11$. Then
\begin{enumerate}
\item[\rm(1)] $M_q= (\sq+1)(q-\sq-2)/4;$
\item[\rm(2)] $M_q'= 3(\sq+1)/2$.
\end{enumerate}
\end{theorem}
\begin{proof} By Lemma \ref{lemma3.2}, 
\begin{equation}\label{eq3.3}
M_q+M_q'=\# \cX(\fq).
\end{equation}
 From
this result, Lemma \ref{lemma3.1} and \S\ref{s2}(\ref{1}),
\begin{equation}\label{eq3.4}
\deg(R_1)=3(2g-2)+\frac{3(\sq+1)}{2}=\frac{\sq-3}{2} M_q'\, .
\end{equation}
The result now follows from (\ref{eq3.3}) and (\ref{eq3.4}), by 
taking into consideration
the maximality of $\cX$ and that $2g-2=(\sq-5)(\sq+1)/4$.
\end{proof}
\section{The example}\label{s4}
In this section we study an example of a plane maximal
curve of degree $\ha(\sq+1)$. In the next section we will see that
this example is, up to $\fq$-isomorphism, the unique plane maximal
curve of degree $\ha(\sq+1)$.

Let $q$ be a square power of a prime $p\ge 3$, and let $\cF$ be the Fermat
curve given by (\ref{eq1.2}). Then $\cF$ is non-singular and maximal.
This is because $\cF$ is covered
by the Hermitian curve  with equation
 $$u^{\sq+1}+v^{\sq+1}+w^{\sq+1}=0$$
 via the
morphism $(u,v,w)\mapsto (U,V,W)=(u^2,v^2,w^2)$ \cite[Prop. 6]{lachaud}.
\begin{remark}\label{rem4.1} (i) The non-inflexion points of $\cF$
relative to $\Sigma_1$ are the ones over $U=\lambda$, over $V=\lambda$
and over $W=\lambda$ for $\lambda$ a $(\sq+1)/2$-th root of $-1$. To see this
we observe that the morphism $U:\cF\to \P^1(\bar\fq)$ has $(\sq+1)/2$
points, say $Q_1,\ldots,Q_{(\sq+1)/2}$ over $U=\infty$ and it
has just one point, say $P_i$, over $U=\lambda_i$ with
$\lambda_i^{(\sq+1)/2}=-1$. Hence, for each
$i=1,\ldots,(\sq+1)/2$, $\div(U-\alpha_i)=\frac{\sq+1}{2}P_i-
\sum_j Q_j$. A similar result holds for $\div(V-\alpha_i)$ and
$\div(W-\alpha_i)$.

(ii) We then see that the Weierstrass semigroup at any of the $3(\sq+1)/2$
points above is $\langle \frac{\sq-1}{2}, \frac{\sq+1}{2}\rangle$. Since
this semigroup cannot be the Weierstrass semigroup at a point of the
non-singular model $\cX$ of $y^{\sq}+y=x^{(\sq+1)/4}$, $\sq\equiv
3\pmod{4}$,
\cite{g-vi}, we conclude that $\cF$ is not $\bar\fq$-isomorphic to $\cX$;
hence these curves are not $\fq$-isomorphic.
\end{remark}

Let $\lambda_1, \ldots, \lambda_{(\sq-1)/2},\lambda:=
\lambda_{(\sq+1)/2}$ be
the roots of $T^{(\sq+1)/2}=-1$, and so each
$\lambda_i$ is in
$\fq$. Let $\cY$ be the non-singular model of the affine curve with equation
\begin{equation}\label{eq4.2}
X^{(\sq+1)/2}= F(Y)\, ,
\end{equation}
with $F(Y)\in \fq [Y]$ satisfying the following properties:
\begin{enumerate}
\item[(a)] $ \deg F =(\sq-1)/2$;
\item[(b)] the  roots of $F$ are $c_j:=(\lambda_j-\lambda)^{-1}$, $j=1,\ldots,(\sq-1)/2$;
\item[(c)] either $F(0)^{\sq-1}=1$ or $F(0)^{\sq-1}=-1$.
\end{enumerate}
\begin{proposition}\label{prop4.1} The curve $\cF$ is $\fq$-isomorphic to $\cY$.
\end{proposition}
\begin{proof}
Write $f= U^{(\sq+1)/2}=\sum_{j=0}^{(\sq+1)/2}A_j(U-\lambda)^j$
with
$A_j=(D^j_U f)(\lambda)$ and  $D^j_U$  the $j$-th Hasse derivative. We
have that $A_0=-1$ and $A_{(\sq+1)/2}=1$, so that
\begin{equation}\label{eq4.3}
\frac{U^{(\sq+1)/2}+1}{{(U-\lambda)}^{(\sq+1)/2}}
=\sum_{j=1}^{(\sq+1)/2}A_j\frac{1}{(U-\lambda)^{(\sq+1)/2-j}}\, .
\end{equation}
Also, equation (\ref{eq1.2}) with $W=1$ is equivalent to
$$
\left[{\frac{V}{U-\lambda}}\right]^{(\sq+1)/2}=\sum_{j=1}^{(\sq+1)/2}
\frac{-A_j}{{U-\lambda}^{(\sq+1)/2-j}}\, .
$$
Consequently for $X=V/(U-\lambda)$ and $Y=1/(U-\lambda)$ we obtain an
equation of type (\ref{eq4.2}). {}From (\ref{eq4.3}),
$$
F(Y)=\sum_{j=1}^{(\sq+1)/2}(-A_j)=
-Y^{(\sq+1)/2}\left[\left(\frac{1}{Y}+\lambda\right)^{(\sq+1)/2}+1\right]
$$
belongs to $\fq [Y]$, it has degree $(\sq-1)/2$ , its roots are
$(\lambda_j-\lambda)^{-1}$
($j=1,\ldots,(\sq-1)/2$), and $F(0)=A_{(\sq+1)/2}\in {\mathbb F}_{\sq}$.

Conversely, let us start with (\ref{eq4.2}). Writing
$F(Y)=k\prod_{j=1}^{(\sq-1)/2}(Y-c_j)$ with $k\in \fq^*$,
$c_j:=\lambda_j-\lambda$, and setting $X=V/(U-\lambda)$ and
$Y=1/(U-\lambda)$, from (\ref{eq4.2}) we find that
$$
V^{(\sq+1)/2}= k(-1)^{(\sq-1)/2}\prod_j c_j(U^{(\sq+1)/2}+1)\, .
$$
Since $k(-1)^{(\sq-1)/2}\prod_j c_j = F(0)=:c^{-1}$, we then have an
equation of type
\begin{equation}\label{eq4.4}
cV^{(\sq+1)/2}=U^{(\sq+1)/2}+1 \mbox{  with } c^{2(\sq-1)}=1.
\end{equation} 
Let
$\epsilon\in \bar\mathbb F_p$ such that $c\epsilon^{(\sq+1)/2}=-1$. Then
(\ref{eq4.4}) implies that $\epsilon\in \fq^*$. Then setting $V=\epsilon V'$ we obtain
an equation of type (\ref{eq1.2}) with $W=1$.
\end{proof}
\section{Proof of the main result}\label{s5}
Throughout the whole section we let $q\ge 121$ and fix the following
notation: 
\begin{enumerate}
\item[(a)] $\cX$ is a non-singular plane maximal curve over $\fq$  
of degree $\ha(\sq+1)$;
\item[(b)] $f=0$ is a minimal equation of $\cX$ with $f\in{\fq}[X,Y]$.
\end{enumerate}

From Lemma \ref{lemma3.1} and Corollary \ref{cor3.2}, $\cX$ has the
following properties: 

\begin{enumerate}
\item[\rm (i)] $\cX$\ is classical for $\Sigma_1;$
\item[\rm (ii)] $\cX$\ is non-classical for $\Sigma_2;$
\item[\rm (iii)] $\cX$\ is Frobenius non-classical for $\Sigma_2.$
\end{enumerate}
Plane curves satisfying (i), (ii), (iii)
above have been characterized 
in terms of their equations \cite{garcia-voloch}, \cite{hk1}.
     \begin{lemma}\label{lemma5.1}
There exist $h, s, z_0,\ldots,z_5\in \fq[X,Y]$ such that
\begin{align}
hf & = z_0^{\sqrt q}+z_1^{\sqrt q}X+
z_2^{\sqrt q}Y+z_3^{\sqrt q}X^2+z_4^{\sqrt q}XY+
z_5^{\sqrt q}Y^2\label{eq5.1}\\
\intertext{and}
sf & = z_0+z_1X^{\sqrt q}+z_2Y^{\sqrt q}+z_3X^{2\sqrt q}+
z_4(XY)^{\sqrt q}+z_5Y^{2\sqrt q}\, .\label{eq5.2}
\end{align}
For a point $P=(a,b,1)\in\cX$ such that $z_i(a,b)\neq 0$ for at
least one index $i,\ 0\leq i\leq 5,$ the conic with equation
$$
z_0(a,b)+z_1(a,b)X+z_2(a,b)Y+z_3(a,b)X^2+z_4(a,b)XY+z_5(a,b)Y^2=0
$$
is the osculating conic of $\cX$ at $P$.
        \end{lemma}
Note that equation (\ref{eq5.2}) is invariant under any change of
projective coordinates. 
To see how the polynomials $z_i$ change, we introduce the matrix
\begin{equation}\label{eq5.3}
\Delta(z_0,\ldots,z_5)=
\left(\begin{array}{ccc}
2z_0 & z_1 & z_2\\
z_2 & 2z_3 & z_4\\
z_3 & z_4 & 2z_5
\end{array}\right)\, ,
\end{equation}
and use homogeneous coordinates $(X)=(X_0,X_1,X_2)$. Now, if the 
change from $(X)$ to $(X')$ is given by $(X)=A(X')$ where $A$ is
a non-singular matrix over $\bar\fq$, then (\ref{eq5.2}) becomes,
again in non-homogeneous coordinates,
\begin{equation}\label{eq5.4}
HF=Z_0^{\sqrt q}+Z_1^{\sqrt q}X'+Z_2^{\sqrt q}Y'+Z_3^{\sqrt q}{X'}^2+
Z_4^{\sqrt q}X'Y'+
Z_5^{\sqrt q}{Y'}^2,
\end{equation}
where $H,F,Z_0,\ldots,Z_5\in {\bar\fq}[X',Y']$ and
$F=0$ is the equation of $\cX$ with respect to the new coordinate
system. Also,
\begin{equation}\label{eq5.5}
\Delta(Z_0,\ldots Z_5)=B^{tr}\Delta(z_0,\ldots,z_5)B\, ,
\end{equation}
where $B$ is the matrix satisfying $B^{\sqrt q}=A$.
If $A$ is a matrix over $\fq$, then $Z_0,\ldots,Z_5\in\fq[X',Y']$,
and (\ref{eq5.1}) becomes
\begin{equation}\label{eq5.6}
SF=Z_0+Z_1{X'}^{\sqrt q}+Z_2{Y'}^{\sqrt q}+Z_3{X'}^{2\sqrt q}+
Z_4(X'Y')^{\sqrt q}+Z_5{Y'}^{2\sqrt q}\, .
\end{equation}

For a rational function $u\in \bar\fq(\cX)$, the symbol $v_P(u)$ denotes the
order of $u$ at $P\in\cX$. Note that $z_i$, for $0\leq i\leq 5$,
can be viewed as a rational function of $\bar\fq(\cX)$. We 
define $e_P:=-{\rm min}_{0\leq i\leq 5} v_P(z_i)$.
        \begin{lemma}\label{lemma5.2}
For $P\in\cX,$ the order $v_P(\det(\Delta(z_0,\ldots,z_5)))$ is 
either $2+e_P$ or $e_P$ according as $P$ is an inflexion point or not.
        \end{lemma}
        \begin{proof}
Take $P$ as the origin and the tangent to $\cX$ at $P$ as the
$X$-axis. Since $P$ is a non-singular point of $\cX$, there
exists a formal power series $y(x)\in{\bar\fq}[[x]]$ of order
$\geq 1$, such that $f(x,y(x))=0$. {}For $0\leq i\leq 5$, put
$m_i=z_i(x,y(x))x^{e_P}$, so that $v_P(m_i(x))\geq 0$. {}From (\ref{eq5.1})
$$
m_0(x)^{\sqrt q}+m_1(x)^{\sqrt q}x+m_2(x)^{\sqrt q}(y(x))+
m_3(x)^{\sqrt q}x^2+m_4(x)^{\sqrt q}(xy(x))+
m_5(x)^{\sqrt q}y(x)^2=0\, .
$$
Putting $y=c_sx^s+\ldots$, with $c_s\not=0$ and $k_i=v_P(m_i(x))$, the
left-hand side is the sum of six formal power series in
the variable $x$ whose orders are as follows:
$$
k_0\sqrt q,\  k_1\sqrt q+1,\  k_2\sqrt q+s,\  k_3\sqrt q+2,\ k_4\sqrt q+s+1,\ 
k_5\sqrt q+2s.
$$
At least two of these orders are equal, and they are less than or equal
to the remaining four. Because of Lemma \ref{lemma3.2} we have two 
possibilities:
\begin{enumerate}

\item[(1)] $s=\ha(\sq+1)$, that is, $P$ is an inflexion point, and $k_0\geq
2,\  k_1=1,\  k_2\geq 1,\  k_3\geq 1,\  k_4\geq 1,\  k_5=0$;

\item[(2)] $s=2$, that is, $P$ is a regular point, and $k_0\geq 1,\ 
k_1\geq 1,\  k_2=k_3=0,\  k_4\geq 0,\  k_5\geq 0$.

\end{enumerate}
In case (1),
$\det(\Delta(z_0(x),\ldots,z_5(x)))=x^{e_P}[cx^2+\ldots]$,
where $c=-c_5c_1^2$ with $m_5(x)=c_5+\ldots$ and $m_1(x)=c_1x+\ldots$.
In  case (2),
$\det(\Delta(z_0(x),\ldots,z_5(x)))=
x^{e_P}[c+\ldots]$, where $c=-c_3c_4$ with $m_3(x)=c_3+\ldots$, and
$m_4(x)=c_4+\ldots$. This completes the proof of the lemma.
        \end{proof}

Following \cite[\S1]{sv}, let $\phi:\cX\to\P^5(\bar\fq)$ be 
the morphism where $\phi(Q)=(z_0,\ldots,z_5)$, for a point $Q\in\cX$,
 and $z_i\in\bar\fq(\cX)$. 
%
%
%
%
Since $P\in\cX$ is a non-singular point of $\cX$, there exists a formal
power series $y(x)\in{\bar\fq}[[x]]$ of order $\geq 1$ such
that $f(x+a,y(x)+b)=0$, where $P=(a,b,1)$. Let
$$
m_i(x)=z_i(x+a,y(x)+b)x^{e_P}\, ,
$$
with $i=0,\ldots,5$. Then we have
\begin{equation*}
\phi(P)=(m_0(x),\ldots,m_5(x))\, ,
\end{equation*}
which is a primitive branch representation of $\phi(P)$.
        \begin{lemma}\label{lemma5.3} The degree of
 $\phi(\cX)$  is $\sq+1$.   
        \end{lemma}
        \begin{proof}
Let $\Sigma$ denote the cubic hypersurface in $\P^5(\bar\fq)$ given by
(\ref{eq5.3}). 
%
%
By the previous lemma, the intersection multiplicity $I(\phi(\cX),\Sigma;
\phi(P))$ of $\phi(\cX)$ and $\Sigma$ at $\phi(P)$ is either $2$ or $0$
according as $P$ is an inflexion point or a regular point of $\cX$. This
shows that $\phi(\cX)$ is not contained in $\Sigma$. {}From
B\'ezout's theorem and Theorem \ref{thm3.2}(2), we obtain 
$3\,\deg(\phi(\cX))=2.3(\sq+1)/2$, whence $\deg(\phi(\cX))=\sq+1$.
        \end{proof}
        \begin{lemma}\label{lemma5.4}
For a generic point $P\in\cX,$ there exists a hyperplane $H$ such that
\begin{enumerate}
\item[(1)] $I(\phi(\cX),H;\phi(P))\geq\sqrt q;$
\item[(2)] the Frobenius image $\Phi(\phi(P))$ lies on $H$.
\end{enumerate}
        \end{lemma}
        \begin{proof}
Choose a point $P=(a,b,1)\in\cX$ such that $z_i(a,b)\neq 0$ for at 
least one index $i$, with $0\le i\le 5$. Then $\phi(P)=(z_0(a,b),
z_1(a,b), z_2(a,b), z_3(a,b), z_4(a,b), z_5(a,b))$. Note that
all points of $\cX$, apart from a finite number of them, are
of this kind. Let $H$ be the hyperplane with equation
$X_0+\alpha X_1+\beta X_2+\alpha^2 X_3+\alpha\beta X_4+\beta^2 X_5=0$,
where $\alpha=a^{\sqrt q}$ and $\beta=b^{\sqrt q}$. There exists a formal
power series $y(x)$ of order $\ge 1$ such that
$f(x+a,y(x)+b)=0$. Putting $z_i(x)=z_i(x+a,y(x)+b))$, we have
$$I(\phi(\cX),\Sigma;\phi(P))= \ord\{z_0(x)+\alpha z_1(x)+\beta z_2(x)+
{\alpha}^2 z_3(x)+{\alpha}{\beta}z_4(x)+{\beta}^2 z_5(x)\}.$$ 
{}From (\ref{eq5.2}) we have
\begin{align*}
& z_0(x)+z_1(x){(x+a)}^{\sqrt q}+z_2(x){(y(x)+b)}^{\sqrt q}+
z_3(x){(x+a)}^{2\sqrt q}+\\
& +z_4(x){((x+a)(y(x)+b))}^{\sqrt q}+z_5(x){(y(x)+b)}^{2\sqrt q}=0\, .
\end{align*}
Since $y(x)$ has order $\ge 1$, that is, $y(x)=cx+\ldots$, then
$$z_0(x)+z_1(x)a^{\sqrt q}+z_2(x)b^{\sqrt q}+
z_3(x)a^{2\sqrt q}+z_4(x)(ab)^{\sqrt q}+z_5(x)b^{2\sqrt q}+
x^{\sqrt q}[\ldots]=0,$$ which proves (1). 

To check (2), note that (\ref{eq5.1}) yields 
$$
z_0(a,b)^{\sqrt q}+z_1(a,b)^{\sqrt q}a+z_2(a,b)^{\sqrt q}b+
z_3(a,b)^{\sqrt q}a^2+z_4(a,b)^{\sqrt q}(ab)+
z_5(a,b)^{\sqrt q}b^2=0\, .
$$
Thus
$$
z_0(a,b)^q+z_1(a,b)^qa^{\sqrt q}+z_2(a,b)^qb^{\sqrt q}+
z_3(a,b)^qa^{2\sqrt q}+z_4(a,b)^q(ab)^{\sqrt q}+
z_5(a,b)^qb^{2\sqrt q}=0\, .
$$
Since $$\Phi(\phi(P))=(z_0(a,b)^q,z_1(a,b)^q,z_2(a,b)^q,z_3(a,b)^q,
z_4(a,b)^q,z_5(a,b)^q),$$ and $\alpha=a^{\sqrt q}, \beta=b^{\sqrt q}$, so
(2) follows.
            \end{proof}
%
%

Now, the linear series of hyperplanes sections of $\phi(\cX)$ is equivalent
to the base-point-free linear series $\cD-E$, where $\cD\cong \P(\langle
z_0,\ldots,z_5\rangle)$ and $E:=\sum_{P\in\cX}e_P P$. By Lemma
\ref{lemma5.3}, this linear series is contained in
$\cD_{\cX} = |(\sq+1)P_0|$,
$P_0\in \cX(\fq)$, because $\cX$ is maximal; hence $(\sq+1)P_0\sim \sq
P+\frx(P)$ (\cite[Corollary 1.2]{fgt}). Note that we do not assert that
equality holds. In fact, this is the case if and only if $\phi(\cX)$ is not
degenerate, that is, $z_0,\ldots, z_5$ are $\bar\fq$-linearly independent.
This gives the following result.
        \begin{lemma}\label{lemma5.5}
The base-point-free linear series of $\cX$ generated by 
$z_0,\ldots,z_5$ is contained in $\cD_{\cX}$.
        \end{lemma}
The next step is to determine the degrees of the $z_i$.
        \begin{lemma}\label{lemma5.6} The degrees satisfy 
${\rm max}_{0\le i\le 5}\deg(z_i)=2$. 
        \end{lemma}
        \begin{proof}
We have seen that the base-point-free linear series 
$\sum_{i=0}^5c_iz_i-E$ on $\cX$ is contained in $\cD_{\cX}$; hence it is
contained in the linear series cut out by
conics on $\cX$, by Theorem \ref{thm3.1}. This implies the existence of
constants 
${d_j}^{(i)}$ such that $\div(z_i)-E= \div(d_i)$, $i=0,\ldots,5$, 
where
$$
d_i=d_i(X,Y)={d_0}^{(i)}+{d_1}^{(i)}X+{d_2}^{(i)}Y+
{d_3}^{(i)}X^2+{d_4}^{(i)}XY+{d_5}^{(i)}Y^2\, .
$$
Choose an index $k$ such that $z_k(X,Y)\not\equiv 0\pmod{f(X,Y)}$. 
Then $\div(z_i/z_k)=\div(d_i/d_k)$. Thus 
$z_i(X,Y)d_k(X,Y)\equiv z_k(X,Y)d_i(X,Y)\pmod{f(X,Y)}$. Now,
re-write (\ref{eq5.1}) in terms of $d_i(X,Y)$:
$$
hfd_k={z_k}^{\sqrt q}({d_0}^{\sqrt q}+d_1^{\sqrt q}X+d_2^{\sqrt q}Y+
d_3^{\sqrt q}X^2+d_4^{\sqrt q}XY+d_5^{\sqrt q}Y^2)\, .
$$
Since $z_k(X,Y)\not\equiv 0\pmod{f(X,Y)}$, so $f(X,Y)$ must divide the other 
factor on the right-hand side, and hence there exists $g\in{\bar\fq}[X,Y]$ such 
that 
$$
gf={d_0}^{\sqrt q}+d_1^{\sqrt q}X+d_2^{\sqrt q}Y+
d_3^{\sqrt q}X^2+d_4^{\sqrt q}XY+d_5^{\sqrt q}Y^2\, ,
$$
with $\deg(d_i)\le 2$, for $i=0,\ldots 5$. Thus we may assume that $g=h$ and
$d_i(X,Y)=z_i(X,Y)$ all $i$. It remains to show that
at least one of the polynomials $z_i(X,Y)$ has degree $2$.
However, if $\deg(z_i(X,Y))\le 1$ for all $i$, then the linear series generated
by $z_0,\ldots, z_5$ would be contained in the linear series cut out by lines.
 But this would imply that $\deg(\phi(\cX))\le (\sq+1)/2$, 
contradicting Lemma \ref{lemma5.3}.
        \end{proof}
%
        \begin{lemma}\label{lemma5.7} 
The polynomials $h$ and $s$ in Lemma \ref{lemma5.1} may be assumed to be
equal.
        \end{lemma}
        \begin{proof}
Since $\deg(z_i)\le 2$ for all $i$, we can re-write
$$
{z_0}+z_1X^{\sqrt q}+z_2Y^{\sqrt q}+
z_3X^{2\sqrt q}+z_4(XY)^{\sqrt q}+z_5Y^{2\sqrt q}
$$ 
in the form
$$
w_0^{\sqrt q}+w_1^{\sqrt q}X+w_2^{\sqrt q}Y+
w_3^{\sqrt q}X^2+w_4^{\sqrt q}XY+w_5^{\sqrt q}Y^2\, ,
$$
where $w_i\in\fq[X,Y]$ and ${\rm max}_{0\le i\le 5}\deg(w_i)= 
{\rm max}_{0\le i\le 5}\deg(z_i)$. Comparing this with (\ref{eq5.1}) 
we see that $z_i$ and $w_i$ only differ 
by a constant in $\fq$ independent of $i$, $0\le i\le 5$. 
Substituting  $cz_i$ for $w_i$  then gives 
\begin{eqnarray}
w_0^{\sqrt q}+w_1^{\sqrt q}X+w_2^{\sqrt q}Y+
w_3^{\sqrt q}X^2+w_4^{\sqrt q}XY+w_5^{\sqrt q}Y^2 & = & \nonumber \\
c^{\sqrt q}({z_0}^{\sqrt q}+z_1^{\sqrt q}X+z_2^{\sqrt q}Y+
z_3^{\sqrt q}X^2+z_4^{\sqrt q}XY+z_5^{\sqrt q}Y^2) & = & c^{\sqrt q}hf \,
.\label{eq5.6a}
\end{eqnarray}
Now, by the previous lemma we
can write $z_i$ explicitly in the form
\begin{equation}\label{eq5.7}
z_i=t_0^{(i)}+t_1^{(i)}X+t_2^{(i)}Y+
t_3^{(i)}X^2+t_4^{(i)}XY+t_5^{(i)}Y^2\, ,
\end{equation}
for $i=0,\ldots 5$. Let $t:=c^{\sq}h$; then (\ref{eq5.6a})  
yields that ${(t_j^{(i)})}^{\sqrt q}=c{t_i}^{(j)}$ for 
$0\le i,j\le 5$. Putting $i=j$, this gives $c^{\sqrt q +1}=1$. Choose 
an element $k$ in $\bar\fq$ such that $k^{\sqrt q-1}=c$, and put 
$d_i=k^{-1}z_i$, $0\le i\le 5$. Then (\ref{eq5.1}) and (\ref{eq5.2}) 
become respectively
\begin{align*}
& hk^{-\sqrt q}f=d_0^{\sqrt q}+d_1^{\sqrt q}X+d_2^{\sqrt q}Y+
d_3^{\sqrt q}X^2+d_4^{\sqrt q}XY+d_5^{\sqrt q}Y^2\ ,\\
\intertext{and}
& tk^{-1}f=k(d_0+d_1X^{\sqrt q}+d_2Y^{\sqrt q}+d_3X^{2\sqrt q}+
d_4(XY)^{\sqrt q}+d_5Y^{2\sqrt q}.
\end{align*}
 Put $h'=hk^{-\sqrt q}$ and $t'=tk^{-1}$. Then
$h'=t'$, and this completes the proof.
        \end{proof}
Next we determine explicitly the coefficients
$t_j^{(i)}$ given in (\ref{eq5.7}), or equivalently the 
$6\times 6$ matrix $T=(t_j^{(i)})$. {}From Lemma \ref{lemma5.7} we
can assume that
\begin{equation}\label{eq5.8}
{(t^{(i)}_j)}^{\sq}=t^{(j)}_i\, .
\end{equation}
for $0\le i,j\le 5$. In other words, we can assume that $T$ is a Hermitian 
matrix over ${\mathbb F}_{\sqrt q}$. 

To obtain further relations between 
elements of $T$, we go back to (\ref{eq5.3}) and note that
$$
{(\det(\Delta(z_0,\ldots,z_5)))}^{\sqrt q}=0
$$
can actually be regarded as the equation of the Hessian curve
$\cH(Z)$ associated to the algebraic curve $\cZ$ with equation
$$z^{\sq}_0+z^{\sq}_1X+{z_2}^{\sqrt q}Y+
{z_3}^{\sqrt q}X^2+{z_4}^{\sqrt q}XY+{z_5}^{\sqrt q}Y^2=0;$$ here
$z_i=z_i(X,Y)$. 
Hence $\cH(\cZ)$ is $\sqrt q$-fold covered by the curve
$\cC$ with equation $\det(\Delta(z_0,\ldots,z_5))=0$, and Lemma
\ref{lemma5.2} 
can be interpreted in terms of intersection multiplicities 
between $\cC$ and $\cX$; namely, $I(\cC,\cX;P)$ is either 
$2+e_P$ or $e_P$ according as $P\in\cX$ is an inflexion 
point or not. Now,
$I(\cH(\cX),\cX;P)=s(P)-2$, where $\cH(\cX)$ is the 
Hessian of $\cX$ and $s(P):= I(\cX,l;P)$, with $l$  
the tangent to $\cX$ at the point $P$; see, for example, \cite[Ch.4,
\S6]{wa}) and, for a characteristic-free
approach to Hessian curves, see \cite[Ch.17]{OO}).  
Comparing the intersection divisors $\cC.\cX$ and 
$\cH(\cX).\cX$, we see that $\frac{n-2}{2}\cC.\cX\ge 
\cH(\cX).\cX$ with $n=\ha(\sq+1)$. Hence, by Noether's ``$AF+BG$'' Theorem,
\cite[p. 133]{Sei}, we obtain
$$
{(\det(\Delta(z_0,\ldots,z_5)))}^{(n-2)/2}=AF+BG\ ,
$$
with $F$ the projectivization of $f$ and $A$, $B$, $G$ 
homogeneous polynomials in ${\bar\fq}[X_0,X_1,X_2]$, 
where $G=0$ is the equation of $\cH(\cX)$. As $\det(
\Delta(z_0,\ldots,z_5))$ is 
a polynomial of degree $6$ (cf. Lemma \ref{lemma5.6}), while
$\deg(G)=3(n-2)$, so $B$ must be a
constant. This yields that $e_P=0$ for each $P\in\cX$. {}For an inflexion
point
$P\in\cX$, we can now infer from the proof of Lemma \ref{lemma5.2} that if 
$P=(0,0,1)$ and $l$ is the $X$--axis, then 
$z_i(0,0)=0$, $i=0,\ldots 4$, and thus $\det(\Delta(z_0,\ldots,z_5))$ has  
no terms of degree $\le 2$. This shows that each inflexion point $P$
of $\cX$ is a singular point of $\cC$. 

By a standard argument depending
on the upper bound $(m-1)(m-2)/2$ for the number of
singular points of an absolutely irreducible algebraic curve of degree
$m$, it can be shown that $\cC$ is doubly covered by an 
absolutely irreducible cubic curve $\cU$ of equation 
$u=0$, with $u$ homogeneous in ${\bar\fq}[X_0,X_1,X_2]$. Hence,
\begin{equation}\label{eq5.9}
\det(\Delta(z_0,\ldots,z_5))=u^2\, .
\end{equation}

Consider now a minor $\Delta_{ij}$ of $\Delta(z_0,\ldots,z_5)$, and
suppose that $\Delta_{ij}$ is not the zero polynomial. 
Then $\Delta_{ij}=0$ can be regarded as the equation of a quartic
curve ${\cV}_{ij}$. Since ${\cV}_{ij}$ also passes through
each inflexion point of $\cX$, so ${\cV}_{ij}$ and $\cU$ have 
at least $3n$ common points. On the other hand, $\deg({\cV}_{ij})
\deg(\cU)=12$, and because $3n>12$, so $\cU$ is a component of 
${\cV}_{ij}$. This shows the existence of linear homogeneous
polynomials $l_0,\ldots,l_5\in {\bar\fq}[X_0,X_1,X_2]$ such that
\begin{alignat}{3}
4z_3z_5-z^2_4 & = ul_0\, ,&
\qquad 
2z_1z_5-z_2z_4 & =-ul_1\, ,&
\qquad
z_1z_4-2z_2z_3 & =ul_2\, ,\label{eq5.10} \\
4z_0z_5-z^2_2 & =ul_3\, ,&
\qquad
2z_0z_4-z_1z_2 & =-ul_4\, ,&
\qquad
4z_0z_3-z^2_1 & =ul_5\, .\label{eq5.11}
\end{alignat}
Let $L$ denote the matrix $\Delta(l_0,l_1,l_2,l_3,l_4,l_5)$. {}From
elementary linear algebra, $\Delta^*=uL$ where $\Delta^*$ is the
adjoint of $\Delta(z_0,\ldots,z_5)$, and hence ${(\det(\Delta(z_0,\ldots,
z_5)))}^2=u^3\det(L)$. Comparison with (\ref{eq5.9}) gives $u=\det(L)$. 
Thus $\Delta^*=\det(L)L$.
%
%
Also, $\Delta(z_0,\ldots,z_5)=\det(L)L^{-1}$; that is,
\begin{alignat}{3}
2z_0 & =l_3l_5-l^2_4\, , &
\qquad
z_1 & =-(l_1l_5-l_2l_4)\, ,&
\qquad
z_2 & =l_1l_4-l_2l_3\, ,\label{eq5.14}\\
2z_3 & =l_0l_5-l^2_2\, ,&
\qquad
z4 & =-(l_0l_4-l_1l_2)\, ,&
\qquad
2z_5 & = l_0l_3-l^2_1\, .\label{eq5.15}
\end{alignat}
Note that we have also seen that $\cU$ has equation $\det(L)=0$. 

Set
\begin{eqnarray*}
 l_i& = & a_iX+b_iY+c_i, \mbox{ for $i=0,2,3,5$,}\\
l_i& = &-a_iX-b_iY-c_i, \mbox{ for $i=1,4$.}
\end{eqnarray*} 
By (\ref{eq5.7}) and (\ref{eq5.14}), we have the following
relations:
\begin{alignat*}{3}
(c_3c_5-c^2_4)/2 & =t_0^{(0)}\, ,&
\qquad
(a_3c_5+c_3a_5-2a_4c_4)/2 & = t_1^{(0)}\, ,&
\quad 
(b_3c_5+c_3b_5-2b_4c_4)/2 & = t_2^{(0)}\, ,\\
(a_3a_5-a^2_4)/2 & =t_3^{(0)}\, ,&
\quad
(a_3b_5+b_3a_5-2a_4b_4)/2 & =t_4^{(0)}\, ,&
\quad 
(b_3b_5-b^2_4)^2 & =t_5^{(0)}\, ;\\
c_1c_5-c_2c_4 & =t_0^{(1)}\, ,&
\quad
a_1c_5+c_1a_5-a_2c_4-c_2a_4 & =t_1^{(1)}\, ,&
\quad 
b_1c_5+c_1b_5-b_2c_4-c_2b_4 & =t_2^{(1)}\, ,\\
a_1a_5-a_2a_4 & =t_3^{(1)}\, ,&
\quad
a_1b_5+b_1a_5-a_2b_4-b_2a_4 & = t_4^{(1)}\, ,&
\quad 
b_1b_5-b_2b_4 & =t_5^{(1)}\, ;\\
c_1c_4-c_2c_3 & = t_0^{(2)}\, ,&
\quad
a_1c_4+c_1a_4-a_2c_3-c_3a_2 & = t_1^{(2)}\, ,&
\quad 
b_1c_4+c_1b_4-b_2c_3-c_2b_3 & = t_2^{(2)}\, ,\\
a_1a_4-a_2a_3 & = t_3^{(2)}\, ,&
\quad
a_1b_4+b_1a_4-a_2b_3-b_2a_3 & = t_4^{(2)}\, ,&
\quad 
b_1b_4-b_2b_3 & = t_5^{(2)}\, ;\\
(c_0c_5-c^2_2)/2 & = t_0^{(3)}\, ,&
\quad
(a_0c_5+c_0a_5-2a_2c_2)/2 & = t_1^{(3)}\, ,&
\quad 
(b_0c_5+c_0b_5-2b_2c_2)/2 & = t_2^{(3)}\, ,\\
(a_0a_5-a^2_2)/2 & = t_3^{(3)}\, ,&
\quad
(a_0b_5+b_0a_5-2a_2b_2)/2 & = t_4^{(3)}\, ,&
\quad 
(b_0b_5-b^2_2)/2 & = t_5^{(3)}\, ;\\
c_0c_4-c_1c_2 & = t_0^{(4)}\, ,&
\quad
a_0c_4+c_0a_4-a_1c_2-c_1a_2 & = t_1^{(4)}\, ,&
\quad 
b_0c_4+c_0b_4-b_1c_2-c_2b_1 & = t_2^{(4)}\, ,\\
a_0a_4-a_1a_2 & = t_3^{(4)}\, ,&
\quad
a_0b_4+b_0a_4-a_1b_2-b_1a_2 & = t_4^{(4)}\, ,&
\quad 
b_0b_4-b_1b_2 & = t_5^{(4)}\, ;\\
(c_0c_3-c^2_1)/2 & = t_0^{(5)}\, ,&
\quad
(a_0c_3+c_0a_3-2a_1c_1)/2 & = t_1^{(5)}\, ,&
\quad 
(b_0c_3+c_0b_3-2b_1c_1)/2 & =t_2^{(5)}\, ,\\
(a_0a_3-a^2_1)/2 & = t_3^{(5)}\, ,&
\quad
(a_0b_3+b_0a_3-2a_1b_1)/2 & = t_4^{(5)}\, ,&
\quad 
(b_0b_3-b^2_1)/2 & = t_5^{(5)}\, .
\end{alignat*}

Now we take an inflexion point $P\in\cX$ to be the origin and the
tangent of $\cX$ at $P$ to be the $X$-axis. Also, since
$I(\cU,\cX;P)=1$, so $P$ is a non-singular point of $\cU$, and
the tangent to $\cU$ at $P$ is not the $X$-axis. We take
this tangent to be the $Y$-axis. We want to show that the $Y$-axis is a
component of $\cU$. As we have seen before, $z_0(0,0)=\ldots z_4(0,0)$,
but $z_5(0,0)\neq 0$. This implies that $t_0^{(0)}=\ldots
=t_4^{(0)}=0$, 
$t_5^{(0)}\neq 0$. By (\ref{eq5.10}) with $u=\det(L)$, $t_0^{(0)}=\ldots
=t_0^{(4)}=0$, 
$t_0^{(5)}\neq 0$. Note that this means, in terms of the $l_i$, that
$c_2=c_4=c_5=0$, $c_1\neq 0$. Putting $k=t_0^{(5)}$, we obtain
\begin{equation}\label{eq5.16}
z_0(X,Y)=kY^2\, .
\end{equation}
Thus the first relation in (\ref{eq5.10}), with $u=\det(L)$, yields
$c_3=0$. Then an easy computation 
shows that $\det(L)=c^2_1(a_5X+b_5Y)+$ terms of degree $\ge 2$, in
$X,Y$. Thus $a_5\neq 0$ but $b_5=0$, since
we have taken the tangent to $\cU$ at $P$ to be the $Y$-axis. Hence
\begin{equation}\label{eq5.17}
l_5=a_5X\, ,\qquad \text{with $a_5\neq 0$}\, .
\end{equation}
Also, $t^{(1)}_2=0$; then (\ref{eq5.8}) gives $t^{(2)}_1=0$, whence
$a_4=0$. Since $k\neq 0$ and $b_5=0$, we obtain
\begin{equation}\label{eq5.18}
l_4=-b_4Y,\quad b_4\neq 0\, .
\end{equation}
The first relation in (\ref{eq5.10}), again with $u=\det(L)$, together with 
(\ref{eq5.16}) and (\ref{eq5.17}) shows that $2kY^2=l_3a_5X-b^2_4Y^2$, whence
\begin{equation}\label{eq5.19}
l_3=0\, .
\end{equation}

With $b=b^{-\sq}_4$, let $A$ be the diagonal matrix with
${\rm diag}(1,1,b^{\sq})$. Change the coordinates from
$(X)$ to $(X')$ by $(X)=A(X')$. {}From the results quoted
after Lemma \ref{lemma5.1}, equations (\ref{eq5.1}) and (\ref{eq5.2}) become 
(\ref{eq5.4}) and (\ref{eq5.6}), where
\begin{alignat*}{2}
H(X',Y') & = h(b^{\sq}X',b^{\sq}Y')\, ,&
\qquad
F(X',Y') & = f(b^{\sq}X',b^{\sq}Y')\, ,\\
Z_0(X',Y') & = z_0(b^{\sq}X',b^{\sq}Y')=b^{2\sq}k{Y'}^2\, ,&
\qquad 
Z_1(X',Y') & = bz_1(b^{\sq}X',b^{\sq}Y')\, ,\\
Z_2(X',Y') & = bz_2(b^{\sq}X',b^{\sq}Y')\, ,&
\qquad 
Z_3(X',Y') & = z_3(b^{\sq}X',b^{\sq}Y')\, ,\\
Z_4(X',Y') & = bz_1(b^{\sq}X',b^{\sq}Y')\, ,&
\qquad 
Z_5(X',Y') & = b^2z_1(b^{\sq}X',b^{\sq}Y')\, ;
\end{alignat*}
by (\ref{eq5.6}),
$\det(\Delta(Z_0,\ldots,Z_5))=b^2\det(\Delta(z_0,\ldots,z_5))$. 
Since $b^{2\sqrt q}=b^{-2}_4$, and $k=t^{(0)}_5={b^2_4}/2$, this
allows us to assume that
\begin{equation}\label{eq5.20}
z_0(X,Y)=-Y^2/2\, ,
\end{equation}
that is, $t^{(0)}_5=-1/2$. Then $b^2_4=1$ and, by (\ref{eq5.8}),
$c^2_1=1$. 
Note that both $b$ and $c$ are in ${\mathbb F}_{\sqrt q}$. 

The next step
is to show that
\begin{align}
b_0b_4-2b_1b_2 & = 0\, , \label{eq5.21}\\
c_0b_4-2c_1b_2 & = 0\, .\label{eq5.22}
\end{align}
Assume first that $b_4=c_4$. As $t^{(1)}_5=-b_2b_4$ and 
$t^{(5)}_1=-a_1c_1$, so (\ref{eq5.8}) yields $a^{\sqrt q}_1=b_2$. 
Similarly, $-b^{\sqrt q}_1=b_1$ follows from $t^{(2)}_5=-b_1b_4$ and 
$t^{(5)}_2=-b_1c_1$, while ${(-a_1b_1)}^{\sqrt q}=b_0b_4-b_1b_2$ 
comes from $t^{(4)}_5=b_0b_4-b_1b_2$ and  
$t^{(5)}_4=-a_1b_1$, again by (\ref{eq5.8}). Putting
together these relations, (\ref{eq5.21}) follows. Since
$t^{(2)}_4=a_1b_4$ and $t^{(4)}_2=b_4c_0-b_2c_1$, so (\ref{eq5.8}) yields 
${(a_1b_4)}^{\sqrt q}=b_4c_0-b_2c_1$. Since
${(a_1b_4)}^{\sqrt q}=b_2b_4=b_2c_1$, so (\ref{eq5.22}) follows. The proof 
for the case $b_4=-c_4$ is similar, and we omit it. 

From the above results
we infer the following. 
        \begin{lemma}\label{lemma5.8}
If $P\in\cX$ is an inflexion$,$ then $\cU$ has a linear component
through $P$.
        \end{lemma}
        \begin{proof} 
We prove that the $Y$-axis is a linear component of $\cU$. 
Equivalently, we can show that $\cX$ is a factor of $\det(L)$. By 
(\ref{eq5.17}) and (\ref{eq5.19}), we must check that $\cX$ divides
$l_0l_4-2l_1l_2$. By (\ref{eq5.17}) and 
$c_2=0$, this occurs if the polynomial 
$(b_0b_4-2b_1b_2)Y^2+(c_0b_4-2c_1b_1)Y$ is identically zero. 
Hence the result is a consequence of (\ref{eq5.21}) 
and (\ref{eq5.22}).
        \end{proof}

It was shown in Theorem
\ref{thm3.2} that $\cX$ has  $3(\sq+1)/2$ inflexion points altogether,
and each one lies on a linear component of $\cU$. 
        \begin{corollary}\label{cor5.1}
The cubic $\cU$ splits into three pairwise distinct lines.
        \end{corollary}
%
%

Let $Q$ be an inflexion point of $\cX$. Take $Q$ as the infinite point
$Y_{\infty}$ and the 
tangent to $\cX$ at $Q$ to be the line at infinity. 
Write (\ref{eq5.1}) in homogeneous coordinates and change coordinates 
from $(X_0,X_1,X_2)$ to $(X_2,X_1,X_0)$. Then (\ref{eq5.1}) becomes
$$
HF=Z^{\sqrt q}_0{X_0}^2+Z^{\sqrt q}_1X_0X_1+Z^{\sqrt q}_2X_0X_2+
Z^{\sqrt q}_3X_1^2+Z^{\sqrt q}_4X_1X_2+Z^{\sqrt q}_5{X^2}_2\, ,
$$
where $Z_i(X_0,X_1,X_2)=z_{5-i}(X_2,X_1,X_0)$, $i=0,\ldots,5$. As $Q$ 
is now the origin and the tangent to $\cX$ at $Q$ is the $X$-axis, so 
(\ref{eq5.16}) gives $Z_0(X_0,X_1,X_2)=cX^2_2$ with 
a non--zero  constant $c$. {}From $Z_0(X_0,X_1,X_2)=z_5(X_2,X_1,X_0)$, 
we then have that $z_5(X_0,X_1,X_2)=cX^2_0$ and hence 
$z_5(X,Y)=c$ in non-homogeneous coordinates. In terms of
$t_i^{(j)}$, this means that $t_i^{(5)}=0$, for  $1\leq i\leq 5$;
hence $a_1=b_1=0$. Taking (\ref{eq5.8}) into account we also have 
$t_5^{(i)}=0$, $1\leq i\leq 5$, which yields
$b_2b_4=b_2b_3=b_0b_5-{b^2}_2=b_0b_4=b_0b_3=0$. Since $b_3=0$ but
$b_4\neq 0$, we get $a_0=b_0=b_1=b_2=0$. Hence $l_0=a_0X$, 
$l_1=c_1$, $l_2=-a_2X$, $l_3=0$, $l_4=-b_4Y$, 
$l_5=a_5X$, where $c^2_1=b^2_4=1$. Note that substituting each 
$l_i$ with $-l_i$, $z_i$ in (\ref{eq5.14}) does not change. This allows us
to put $c_1=1$. Actually, we can also assume $b_4=1$, as the change of
coordinates $X'=X$, $Y'=-Y$ changes $l_4$ into $-l_4$ while preserving
the other $l_i$. Thus $t_1^{(4)}=-a_2, 
t_4^{(1)}=a_2, t_1^{(1)}=a_5, t_4^{(4)}=a_0$, whence
$a_0,a_2,a_5\in {\mathbb F}_{\sq}$  by (\ref{eq5.8}). 
The following has therefore been shown.
        \begin{lemma}\label{lemma5.9} There exist  $a_0,a_2,a_5\in
{\mathbb F}_{\sq}$ such that
\begin{align}
l_0 &  = a_0X\, ,\quad l_1=1\, ,\quad l_2=-a_2X\, ,\quad l_3=0\, ,\quad
l_4=-Y\, ,
\quad l_5=a_5X\, ,\label{eq5.23}\\
z_0 & =-Y^2/2\, ,\quad z_1=-a_5X+a_2XY\, ,\quad z_2=-Y\, ,\label{eq5.24}\\
z_3 & = (a_0a_5-a^2_2)/2X^2\, ,\quad z_4=a_0XY-a_2X\, ,
\quad z_5=-1/2\, .\label{eq5.25}
\end{align}
        \end{lemma}

Now we want to show  that, if $R=(0,\eta)$ is any further inflexion point 
of $\cX$ lying on the $Y$-axis, then the tangent line $r$ to $\cX$ at 
$R$ has equation $Y=\eta$. To do this it is sufficient to check that the
curve $\cZ$ with equation  
$$z^{\sq}_0+z^{\sq}_1X+z^{\sq}_2Y+
z^{\sq}_3X^2+z^{\sq}_4XY+z^{\sq}_5Y^2=0$$
has a cusp at $R$, that is, a double point with only one tangent,
such that the tangent is the horizontal line $Y=\eta$. Applying the
translation 
$X'=X$, $Y'=Y-\eta$, the curve $\cZ$ is transformed into the curve with
equation
$$
-{(\eta^{\sq}+\eta)}^2/2+(\eta^{\sq}+\eta)Y-Y^2/2+ \alpha=0\, ,
$$
where $\alpha$ represents terms of degree at least $3$. 
Since this curve passes through the origin, we have $\eta^{\sq}+\eta=0$. 
Hence the lowest degree term is $-Y^2/2$ and so the origin is a cusp with
tangent line $Y=0$, as required. This gives the following situation.
        \begin{theorem}\label{thm5.1}
There exists a triangle such that the  inflexion points
 of
$\cX$ lie  $\ha(\sq+1)$ on each side$,$ none a vertex$,$ and the 
inflexional tangents pass $\ha(\sq+1)$ through  each vertex$,$ 
none being a side.
        \end{theorem}

We are now in a position to prove the main result, Theorem \ref{thm1},
stated in \S1.

Let $n=(\sq+1)/2$. We choose the  triangle $\cT$ of Theorem \ref{thm5.1}
as triangle of reference, and denote 
the inflexions on the $X$-axis by $(\xi_i,0)$, $i=1\ldots n$, and those 
on the $Y$-axis by $(0,\eta_i)$, $i=1\ldots n$. Also, without loss of 
generality, we may assume that ${\xi_1}^n+1=0$ and 
${\eta_1}^n+1=0$. Write $f(X,Y)$ in the 
form 
$$
f=a_0(X)Y^n+\ldots +a_j(X)Y^{n-j}+\ldots +a_n(X),
$$
 with $a_i(X)$ of degree $i$ in $\fq[X]$.

Since $(\xi_i,0)$ lies on $\cX$, so $a_n(\xi_i)=0$. 
Since the line $x=\xi_i$ is the inflexional tangent at 
$(\xi_i,0)$, so
$$a_0(\xi_i)Y^n+\ldots+a_{n-1}(\xi_i)Y =0$$  has $n$ repeated roots. So
$$a_1(\xi_i) = \ldots = a_{n-1}(\xi_i)= 0.$$ 
Since is true for all $\xi_i$,
$$a_1(X)=\ldots = a_{n-1}(X) = 0.$$ 
Hence
$f(X,Y)=a_0Y^n+a_n(X)$. A similar argument shows that
$f(X,Y)=b_0X^n+b_n(Y)$. Thus $f(X,Y)=a_0X^n+b_0Y^n+c_0$, and
it only remains to compute the coefficients.
Since $f(\xi_1,0)=0$ and ${\xi_1}^n+1=0$, we have $a_0=c_0$.
Similarly, from ${\eta_1}^n+1=0$ we infer $b_0=c_0$. This completes
the proof.
\bigskip

\noindent {\bf Acknowledgments. } The authors were supported respectively
by a NATO grant, G.N.S.A.G.A. Italian C.N.R., and Cnpq-Brazil.
\bigskip

\end{document}